\documentclass[11pt]{article}

\usepackage[utf8]{inputenc}
\usepackage[a4paper]{geometry}
\usepackage{amsmath,mathtools}
\usepackage{amsfonts,amssymb,latexsym,multicol,color}
\usepackage[francais]{babel}
\usepackage[T1]{fontenc}
\usepackage{graphicx}
\newcounter{fig}
\input{cyracc.def}

\newtheorem{theo}{Th\'eor\`eme}

\newtheorem{prop}{Proposition}

\newcommand{\cad}{\text{c'est-\`a-dire }}
\newcommand{\expli}[1]{\quad\text{\footnotesize (#1)}}

\makeatletter
\ifcase\@ptsize

\or

\or

\fi
\makeatother

\newcommand{\e}{{\rm e}}
\newcommand{\epv}{\quad ; \quad}

\newcommand{\implique}{\Rightarrow}
\newcommand{\inv}{{\rm inv}}
\newcommand{\Kl}{{\rm K}}

\newcommand{\ioe}{\leqslant}
\newcommand{\soe}{\geqslant}
\renewcommand{\le}{\leqslant}
\renewcommand{\ge}{\geqslant}
\newcommand{\vers}{\rightarrow}

\newcommand{\ssi}{\Leftrightarrow}

\newcommand{\mes}{{\rm mes}}

\newcommand{\pepv}{\; ; \;}

\newcommand{\eps}{\varepsilon}
\newcommand{\fhi}{\varphi}
\newcommand{\thet}{\vartheta}

\newcommand{\abs}[1]{\left\lvert #1 \right\rvert}

\newcommand{\set}[1]{\left\{ #1 \right\}}

\newcommand{\Fcal}{{\mathcal F}}

\newcommand{\Nat}{{\mathbb N}}
\newcommand{\Int}{{\mathbb Z}}

\newcommand{\Real}{{\mathbb R}}
\newcommand{\Com}{{\mathbb C}}

\newcommand{\vide}{\text{\O}}

\newcommand{\card}{{\rm card}}

\newcommand{\fin}{\hfill$\Box$}
\newcommand{\dem}{\noindent {\bf D\'emonstration\ }}
\newcommand{\fine}{\tag*{\mbox{$\Box$}}}

\providecommand{\bysame}{\leavevmode ---\ }
\providecommand{\og}{``}
\providecommand{\fg}{''}
\providecommand{\smfandname}{et}
\providecommand{\smfedsname}{\'eds.}
\providecommand{\smfedname}{\'ed.}
\providecommand{\smfmastersthesisname}{M\'emoire}
\providecommand{\smfphdthesisname}{Th\`ese}

\newcommand{\virg}{\raisebox{.7mm}{,}}

\title{Démonstration d'une conjecture de Kruyswijk et Meijer sur le plus petit dénominateur\\ des nombres rationnels d'un intervalle}

\author{Michel Balazard et Bruno Martin\footnote{B. Martin est financé par le projet ANR-FWF : FWF: I 4945-N et ANR-20-CE91-0006.}}

\date{}

\begin{document}
\maketitle

\begin{center}
  {\sc Abstract}
\end{center}
\begin{quote}
{\footnotesize The average value of the smallest denominator of a rational number belonging to the interval~$](j-1)/N,j/N]$, where~$j=1,\dots, N$, is proved to be asymptotically equivalent to~$16\pi^{-2}\sqrt{N}$, when $N$ tends to infinity. The result had been conjectured in 1977 by Kruyswijk and Meijer.}
\end{quote}

\begin{center}
  {\sc Keywords}
\end{center}
\begin{quote}
{\footnotesize Farey sequences, Kloosterman sum }
\end{quote}


\section{Introduction} 

Si $E$ est une partie de $\Real$, notons $q(E)$ le plus petit dénominateur d'une fraction irréductible représentant un nombre rationnel appartenant à $E$ (avec $q(E)=\infty$ si $E$ ne contient aucun nombre rationnel) :
\[
q(E)=\min \set{q \in \Nat^* \pepv \exists\, p \in \Int, \; p/q \in E}.
\]
Nous dirons que $q(E)$ est le \emph{dénominateur minimal} de $E$. Observons que $E \mapsto q(E)$ est une fonction décroissante.

Pour $t \in \Real$ et $\delta >0$, posons
\[
q(t,\delta)=q(\,]t-\delta,t]), 
\]
et pour $N \in \Nat^*$ et $j=1,\dots,N$, posons
\[
q_j(N)=q\Big( \frac jN, \frac 1 N\Big)=q\Big(\,\Big]\frac{j-1}{N},\frac{j}N\Big]\,\Big).
\]
On a, par exemple, $q_1(N)=N$ et $q_N(N)=1$. On pose
\[
S(N)=\sum_{j=1}^Nq_j(N).
\]
Kruyswijk et Meijer ont établi en 1977 l'existence de constantes positives $C_1$ et $C_2$ telles que $C_1 \ioe S(N)/N^{3/2} \ioe C_2$ et conjecturé que
\begin{displaymath}
S(N)\sim  \frac{16}{\pi^2} N^{3/2}\quad(N\to +\infty)
\end{displaymath}
(cf. \cite{zbMATH03573955}).
Stewart  a ensuite obtenu un premier encadrement explicite de $S(N)$ (cf. \cite{zbMATH03534552}), qu'il a amélioré en 2013 dans l'article \cite{zbMATH06204646}  : pour $N$ suffisamment grand, on a
\begin{displaymath}
   1.35 N^{3/2} < S(N) < 2.04 N^{3/2}. 
\end{displaymath}

Nous confirmons la conjecture de Kruyswijk et Meijer sous la forme suivante.
\begin{theo}
  Pour $N\ge 2$, on a 
  \begin{equation}\label{eq:resultat-principal}
    S(N) =\frac{16}{\pi^2} N^{3/2}+O(N^{4/3} \ln^2 N). 
  \end{equation}
\end{theo}

Notre démonstration repose sur une estimation récemment obtenue par Chen et Haynes, dans la prépublication \cite{ch23}, pour une version continue du problème : on a
\begin{equation}\label{eq:resultat-Chen-Haynes}
  \int_0^1 q(t,\delta)dt =  \frac{16}{\pi^2} \frac{1}{\delta^{1/2}} +O(\ln^2 \delta) \quad ( 0 < \delta \ioe 1/2).
\end{equation}
En suivant l'approche proposée par Stewart, nous obtenons une formule explicite pour 
la différence 
\begin{displaymath}
  R(N) = S(N) -N\int_0^1 q(t,1/N)dt 
\end{displaymath}
 et établissons l'estimation
 \begin{equation}\label{eq:majoration-R(N)}
   R(N) =O(N^{4/3}  \ln^2 N),
 \end{equation}
en utilisant, notamment, la majoration de Weil de la valeur absolue de la somme de Kloosterman. Le théorème est alors une conséquence directe de  \eqref{eq:resultat-Chen-Haynes} et \eqref{eq:majoration-R(N)}. 

\smallskip 

Au lieu de choisir l'intervalle semi-ouvert $]t-\delta,t]$ pour la définition de $q(t,\delta)$, on peut considérer l'un des trois autres types d'intervalles possibles : $[t-\delta,t[$, $[t-\delta,t]$ et~$]t-\delta,t[$. Le théorème reste inchangé dans tous les cas. Nous en donnons une justification au paragraphe \ref{sec:variante-th}.

\smallskip

Nous utiliserons parfois la notation d'Iverson : $[P]=1$ si la propriété $P$ est vraie, et~$[P]=0$ si elle est fausse. Rappelons également que la notation de Vinogradov $A\ll B$ signifie que $A=O(B)$. Enfin nous employons la notation $\e(x)=\exp(2i \pi x)$.

\section{Préalables arithmétiques}

\subsection{Fonction plafond}

Pour $x\in \Real$, $\lceil x\rceil$ désigne le plus petit nombre entier supérieur ou égal à $x$. En désignant par~$\set{t}$ la partie fractionnaire du nombre réel $t$, on a $\lceil x\rceil =x + \set{-x}$. Si $a,b \in \Real$, on a
\begin{equation}\label{230201a}
\card ([a,b[ \, \cap \, \Int)=
\begin{dcases}
\quad \;\; 0 & (a > b)\\
\lceil b \rceil - \lceil a \rceil & (a \ioe b).
\end{dcases}
\end{equation}

\subsection{Fonction inverse modulaire}

Si $q$ est un nombre entier positif et $a$ un nombre entier premier à $q$, nous noterons~$\inv(a,q)$ le nombre entier $b$ tel que
\[
0<b \ioe q \; \text{ et } \; ab \equiv 1 \pmod q.
\]

Si le nombre $q$ est clairement déterminé par le contexte, on note souvent $\bar{a}$ au lieu de~$\inv(a,q)$.

\subsection{Rappels sur les ensembles de Farey}\label{230131d}

Pour $k\in \Nat^*$, nous noterons $\Fcal_k$ l'ensemble de Farey d'ordre $k$, constitué des nombres rationnels représentés par des fractions irréductibles $p/q$ avec $p \in \Int$, $q \in \Nat^*$, $q \ioe k$, et~$(p,q)=1$ (dites {\og fractions de Farey d'ordre $k$\fg}), muni de la relation d'ordre usuelle. 

La pertinence des ensembles de Farey pour l'étude du dénominateur minimal d'une partie de~$\Real$ repose sur l'équivalence suivante,
\begin{equation}\label{230131b}
q(E) > k \ssi E\, \cap \, \Fcal_k= \vide.
\end{equation}

\smallskip

Soit 
\begin{equation}\label{230221a}
0= \rho_0 < \rho_1 < \dots < \rho_M=1
\end{equation}
les éléments de $\Fcal_k \, \cap \, [0,1]$.  On a $M=A(k)$, où, en notant $\fhi$ l'indicatrice d'Euler,
\[
A(x)=\sum_{n \ioe x} \fhi(n)  \quad(x >0).
\]

Nous utiliserons les propriétés suivantes. 

\medskip

Si deux fractions de Farey d'ordre $k$, disons $a/r$ et $b/s$ avec $a/r <b/s$, sont consécutives dans $\Fcal_k$, alors
\[
br-as=1 \quad\text{ et }\quad r+s >k.
\]
(Theorems 28 et 30, p. 28 de \cite{zbMATH05309455}).  
Si~$0<b/s \ioe 1$, il résulte de la première de ces relations que~$b=\inv(r,s)$.

Par conséquent, chaque élément $\rho_i$ de \eqref{230221a}, avec $1 \ioe i \ioe M$, est déterminé par le couple, disons $(r,s)$, où $r$ est le dénominateur de $\rho_{i-1}$ et $s$ celui de $\rho_i$. Lorsque $i$ varie de~$i=1$ à~$i=M$, la suite de ces couples décrit bijectivement l'ensemble
\begin{equation}\label{230221b}
\set{(r,s) \in \Nat^{*2} \pepv \max(r,s) \ioe k < r+s, \; (r,s)=1}.
\end{equation}
On observe que 
\begin{displaymath}
  \rho_{i}-\rho_{i-1}= \frac 1{rs}\cdotp
\end{displaymath}

Enfin, si $r$ et $s$ sont les dénominateurs respectifs de deux fractions consécutives, disons~$\rho<\rho'$ dans $\Fcal_k$, alors le dénominateur $t$ de la fraction $\rho''$ suivant immédiatement~$\rho'$ est
\begin{equation}\label{230209b}
t=\thet(k,r,s)=s\Big\lfloor \frac{k+r}s\Big\rfloor -r=k-s\set{\frac{k+r}s}
\end{equation}
(cf. \cite{zbMATH03880795}, Lemma 1, p. 399, et (2), p. 400).

\section{Une application de la transformation de Fourier discrète}

\subsection{Définition et formule d'inversion}

Nous utiliserons la transformation de Fourier discrète (cf. \cite{zbMATH06404910}, Chapter 2).
Soit $q$ un nombre entier positif, et $f : \Int \vers \Com$ une fonction de période $q$. On pose
\[
\widehat{f}(x)=\sum_{n=1}^{q} f(n)\,\e\Big(-\frac{nx}{q}\Big) \quad (x \in \Int).
\]
La fonction $\widehat{f}$ est la transformée de Fourier (discrète) de $f$, relative à sa période $q$. C'est une fonction définie sur $\Int$, de période $q$.  La transformation de Fourier d'une fonction paire (resp. impaire) est paire (resp. impaire).
Connaître $\widehat{f}$ et $q$ permet de retrouver $f$ par la formule d'inversion de Fourier,
\[
f(n)=\frac 1q\sum_{x=1}^{q} \widehat{f}(x)\,\e\Big(\frac{nx}{q}\Big) \quad (n \in \Int).
\]


\subsection{Sommes de Kloosterman et de Ramanujan}

Pour $a,b \in \Int$ et $q\in \Nat^*$, on définit la somme de Kloosterman,
\[
\Kl(a,b\,;\,q)=\sum_{\substack{n=1\\(n,q)=1}}^q\e\Big( \frac{an+b\bar{n}}{q}\Big),
\]
où $\bar{n}=\inv(n,q)$. La valeur de cette somme est réelle.

Si $b=0$, on note $\Kl(a,0\,;\,q)=c_q(a)$, somme de Ramanujan. C'est une fonction paire de $a$, qui n'est autre que la transformée de Fourier de la fonction indicatrice des nombres premiers à~$q$. On a aussi $\Kl(0,b\,;\,q)=c_q(b)$.

Nous utiliserons l'inégalité de Weil,
\[
\abs{\Kl(a,b\,;\,q)} \ioe (a,b,q)^{1/2}\, \tau(q)\, \sqrt{q},
\]
où $\tau$ est la fonction {\og nombre de diviseurs\fg} (cf. \cite{zbMATH02121181}, Corollary 11.12, p. 280).

\subsection{Sommes faisant intervenir l'inversion modulo $q$}

Soit $q$ un nombre entier positif, et $f : \Int \vers \Com$ une fonction de période $q$. Soit $N$ un nombre entier positif. Posons
\[
g(n)=[(n,q)=1] \cdot f(N\bar{n}) \quad (n \in \Int).
\]
La fonction $g$ est de période $q$. Calculons sa transformée de Fourier :
\begin{align*}
\widehat{g}(x) &= \sum_{n=1}^{q} g(n)\e\Big(-\frac{nx}{q}\Big)=\sum_{\substack{n=1\\(n,q)=1}}^qf(N\bar{n})\e\Big(-\frac{nx}{q}\Big)\\
&=\sum_{\substack{n=1\\(n,q)=1}}^q\e\Big(-\frac{nx}{q}\Big)\Big( \frac 1q\sum_{y=1}^{q} \widehat{f}(y)\e\Big( \frac{N\bar{n}y}q\Big) \Big)\\
&=\frac 1q\sum_{y=1}^{q} \widehat{f}(y) \sum_{\substack{n=1\\(n,q)=1}}^q\e\Big( \frac{-nx+N\bar{n}y}q\Big)
=\frac 1q\sum_{y=1}^{q}  \widehat{f}(y) \Kl(-x,Ny\,; \,q).
\end{align*}
Par inversion de Fourier, on en déduit
\begin{align*}
g(n)&=\frac 1q\sum_{x=1}^{q} \widehat{g}(x)\e\Big(\frac{nx}{q}\Big) \\
&=\frac 1q\sum_{x=1}^{q} \e\Big(\frac{nx}{q}\Big) \Big(\frac 1q\sum_{y=1}^{q}  \widehat{f}(y) \Kl(-x,Ny\,; \,q)\Big)\\
&=\frac{1}{q^2} \sum_{y=1}^{q}  \widehat{f}(y)\sum_{x=1}^{q}\e\Big(\frac{nx}{q}\Big) \Kl(-x,Ny\,; \,q).
\end{align*}

En faisant l'hypothèse supplémentaire que $f$ est impaire, nous allons en déduire une majoration de la somme des valeurs de la fonction $g$ sur un intervalle. Notre démonstration est inspirée de celle de l'inégalité de P\'olya-Vinogradov, telle qu'elle est exposée, par exemple, dans~\cite{MR1790423}, Chapter~23.

Soit $I$ un intervalle fini de $\Int$. On a
\begin{align}
&\sum_{\substack{n \in I\\(n,q)=1}} f(N\bar{n}) =\sum_{n \in I} g(n)\notag\\
&=\frac{1}{q^2} \sum_{n \in I}\; \sum_{y=1}^{q-1}  \widehat{f}(y)\sum_{x=1}^{q} \e\Big(\frac{nx}{q}\Big)  \Kl(-x,Ny\,; \,q) \expli{$\widehat{f}(q)=\widehat{f}(0)=0$ car $\widehat{f}$ est impaire}\notag\\
&=\frac{1}{q^2}\sum_{y=1}^{q-1}  \widehat{f}(y)\sum_{x=1}^{q} \Kl(-x,Ny\,; \,q)\sum_{n \in I} \e\Big(\frac{nx}{q}\Big)\notag\\
&=\frac{\card\, I}{q^2}\sum_{y=1}^{q-1}  \widehat{f}(y)c_q(Ny) + \frac{1}{q^2}\sum_{y=1}^{q-1}  \widehat{f}(y)\sum_{x=1}^{q-1} \Kl(-x,Ny\,; \,q)\sum_{n \in I} \e\Big(\frac{nx}{q}\Big).\label{230216a}
\end{align}

D'une part, comme $\widehat{f}$ est impaire et $c_q$ paire, on a
\[
\sum_{y=1}^{q-1}  \widehat{f}(y)c_q(Ny) =0.
\]
D'autre part, en appliquant la majoration 
\begin{equation}\label{230216b}
\Big \lvert \sum_{n \in I} \e(n\alpha)\Big\rvert \ioe \frac{1}{\sin \pi \alpha} \quad ( 0<\alpha<1),
\end{equation}
et la majoration de Weil de la somme de Kloosterman, on obtient, pour $y \in \set{1,\dots, q-1}$,
\begin{align}
\Big\lvert \sum_{x=1}^{q-1} \Kl(-x,Ny\,; \,q)\sum_{n \in I} \e\Big(\frac{nx}{q}\Big) \Big\rvert & \ioe \tau(q)q^{1/2}\sum_{x=1}^{q-1} \frac{(x,Ny,q)^{1/2}}{\sin (\pi x/q)}\notag\\
&\ioe \tau(q)q^{1/2}\sum_{x=1}^{q-1} \frac{(x,q)^{1/2}}{\sin (\pi x/q)}.\label{230216d}
\end{align}
Or on a
\begin{align}
\sum_{x=1}^{q-1} \frac{(x,q)^{1/2}}{\sin (\pi x/q)} &=\sum_{d \mid q} \sqrt{d} \, \sum_{\substack{x=1\\ (x,q)=d}}^{q-1} \frac{1}{\sin (\pi x/q)} \notag\\
&\ioe \sum_{d \mid q} \sqrt{d} \, \sum_{\substack{x=1\\ d \mid x}}^{q-1} \frac{1}{\sin (\pi x/q)} 
=\sum_{d \mid q} \sqrt{d} \, \sum_{m < q/d} \frac{1}{\sin \big(\pi m/(q/d)\big)} \notag\\
&\ioe\sum_{d \mid q} \sqrt{d} \cdot \frac qd\ln (q/d) \expli{cf. \cite{MR1790423}, p. 136}\notag\\
&=q\beta(q),\label{230216c}
\end{align}
disons, où l'on a noté $\beta$ la fonction arithmétique définie par
\begin{equation}
  \label{eq:def-beta}
  \beta(q)=\sum_{d \mid q} \frac{\ln(q/d)}{\sqrt{d}}\cdotp
\end{equation}

\smallskip

En insérant \eqref{230216d} et \eqref{230216c} dans \eqref{230216a}, on obtient, pour $f$ impaire, la majoration
\begin{equation}\label{230216e}
\Big\lvert  \sum_{\substack{n \in I\\(n,q)=1}} f(N\bar{n})  \Big \rvert \ioe \frac{\tau(q)\beta(q)}{\sqrt{q}}\sum_{y=1}^{q-1}  \big\lvert\widehat{f}(y)\big\rvert.
\end{equation}

\smallskip

Une majoration de la fonction sommatoire de la fonction arithmétique $\tau\beta$ nous sera utile :
\begin{align}
\sum_{n\ioe x} \tau(n) \beta(n)  &=\sum_{n \ioe x} \tau(n) \sum_{d \mid n} \frac{\ln(n/d)}{\sqrt{d}}
=\sum_{d \ioe x} \frac{1}{\sqrt{d}}\sum_{m \ioe x/d}\tau(md) \ln m \notag\\
&\ioe \sum_{d \ioe x} \frac{\tau(d)}{\sqrt{d}}\sum_{m \ioe x/d}\tau(m) \ln m 
\ll \sum_{d \ioe x} \frac{\tau(d)}{\sqrt{d}} (x/d)\ln^2(x/d)\notag\\
& \ll x \ln^2 x \quad (x \soe 1).\label{230222a}
\end{align}

\subsection{Application à la première fonction de Bernoulli}

Nous désignons par  $B_1$ la première fonction de Bernoulli, définie par
\[
B_1(x)=
\begin{dcases}
0 & (x \in \Int)\\
\set{x}-1/2 & ( x \notin \Int).
\end{dcases}
\]

\begin{prop}
Soit $I$ un intervalle de nombres entiers, $q$ et $N$ des nombres entiers positifs. On a 
  \begin{equation}\label{230216f}
    \Big\lvert  \sum_{\substack{n \in I\\(n,q)=1}} B_1\Big( \frac{N\bar{n}}{q}\Big)  \Big \rvert \ioe \frac 12\sqrt{q}\,\tau(q)\beta(q)\ln q, 
  \end{equation}
où $\beta(q)$ a été définie en \eqref{eq:def-beta}. 
\end{prop}

\dem 
Calculons la transformée de Fourier discrète de $f(n)=B_1(n/q)$. 
On a
\begin{align*}
\widehat{f}(x) &= \sum_{n=1}^{q}  B_1\Big( \frac n q \Big)\e\Big(-\frac{nx}{q}\Big)\\
&=\sum_{n=1}^{q-1} \Big( \frac nq -\frac 12\Big)\e\Big(-\frac{nx}{q}\Big)\\
&=\frac 1q \sum_{n=1}^{q-1} n\e\Big(-\frac{nx}{q}\Big)-\frac 12\sum_{n=1}^{q-1} \e\Big(-\frac{nx}{q}\Big).
\end{align*}
Si $q \mid x$, on trouve $\widehat{f}(x) =0$. Si $q \nmid x$, l'identité
\[
\sum_{n=1}^{q-1} nz^n=z\frac{(q-1)z^q-qz^{q-1}+1}{(z-1)^2} \quad (z \neq 1),
\]
fournit l'égalité
\[
\widehat{f}(x) =\frac{\e(x/q)}{1-\e(x/q)}+\frac 12=\frac{1+\e(x/q)}{2\big(1-\e(x/q)\big)}\cdotp
\]

\smallskip

Comme $B_1$ est impaire, et comme
\[
\sum_{y=1}^{q-1}  \big\lvert\widehat{f}(y)\big\rvert \ioe\sum_{y=1}^{q-1} \frac{1}{2\sin \pi y/q}\ioe \frac{q\ln q}2\virg
\]
l'inégalité \eqref{230216e} entraîne \eqref{230216f}.
\fin

\smallskip 

Pour évaluer le terme $R(N)$ dans le paragraphe \ref{sec:esti-RN}, nous aurons recours à l'estimation suivante qui découle directement de \eqref{230216f} par sommation partielle. 
\begin{prop}\label{eval-somme-B1-IPP}
Soit $s\in\Nat^*$, $a,b$ deux nombres réels tels que $a<b$,  et $N\in \Nat$. On a 
\begin{displaymath}
  \Big|\sum_{\substack{a<r\le b\\ (r,s)=1}} (r-a) B_1\Big( \frac{\overline{r}N}s\Big)\Big|
 \le (b-a)\sqrt{s}\,\tau(s)\beta(s)\ln s. 
\end{displaymath}
\end{prop}

\section{Identité pour $S(N)$}\label{230224a}

En explicitant l'approche proposée par Stewart dans \cite{zbMATH03534552} et \cite{zbMATH06204646}, nous allons obtenir une formule explicite pour la différence
\[
R(N)=S(N)-N\int_0^1\, q(t,1/N) \, dt
\]
(cf. Proposition \ref{230217b} ci-dessous).

\subsection{Expression de $R(N)$ au moyen des lois de répartition des dénominateurs minimaux}

Nous transformons d'abord la somme $S(N)$ à l'aide d'une manipulation classique. Posons, pour~$1 \ioe k \ioe N$, 
\[
\alpha_N(k)=\card\set{j \in\set{1,\dots,N} \pepv q_j(N)=k}.
\]
 On a donc
 \begin{equation*}
 \sum_{k=1}^N \alpha_N(k) = N \epv  \sum_{k=1}^N k\alpha_N(k) = S(N).
 \end{equation*}
 Posons ensuite, pour $0 \ioe k \ioe N$,
 \[
\theta_N(k)=\card\set{j \in\set{1,\dots,N}  \pepv q_j(N) >k}.
\]
 Ainsi,
 \begin{equation*}
  \theta_N(0) = N \epv \theta_N(N)=0 \epv \alpha_N(k)=\theta_N(k-1)-\theta_N(k) \quad (1 \ioe k \ioe N).
 \end{equation*} 
Par conséquent,
\begin{equation}\label{eq:identite-SN}
S(N) = \sum_{k=1}^N k \alpha_N(k)
=  \sum_{k=1}^Nk\big( \theta_N(k-1)-\theta_N(k)\big)
=\sum_{k=0}^N\theta_N(k).
\end{equation}

\medskip

En posant maintenant, pour $0 \ioe k \ioe N$,
 \[
\nu_N(k)=\mes\set{t\in [0,1[  \pepv q(t,1/N) >k}
\]
(où $\mes$ désigne la mesure de Lebesgue), une démonstration analogue fournit l'égalité
\[
\int_0^1\, q(t,1/N) \, dt=\sum_{k=0}^N\nu_N(k).
\]
En remarquant que $\theta_N(0)-N\nu_N(0)=N-N=0$, on obtient par différence la relation
\begin{equation}\label{230217a}
R(N)= \sum_{k=1}^N\big(\theta_N(k)-N\nu_N(k)\big).
\end{equation}

\subsection{Sous-intervalles de $]0,1]$ ne contenant aucun terme d'une suite donnée}

Pour donner une expression de la différence $\theta_N(k)-N\nu_N(k)$, nous utiliserons le lemme suivant.
\begin{prop}\label{230130a}
Soit $M,N \in \Nat^*$, $\delta >0$, et $F$ un ensemble de $M+1$ points du segment~$[0,1]$, dont les éléments, notés $\rho_i$, $i=0,\dots, M$, vérifient
\begin{equation*}\label{230131c}
\rho_0=0 < \rho_1 < \dots < \rho_M=1.
\end{equation*}
Soit 
\[
X(\delta)=\set{t \in\,  ]0,1]  \pepv ]t-\delta,t] \, \cap \  F= \vide} .
\]
On a
\begin{equation*}
X(\delta) =\bigcup_{i=1}^M\, [\rho_{i-1}+\delta,\rho_i[
\end{equation*}
et
\begin{equation*}
\Big\{j \in \set{1,\dots,N}  \pepv \Big] \frac{j-1}{N}, \frac jN\Big] \, \cap \, F= \vide\Big\} =\bigcup_{i=1}^M\, ([N\rho_{i-1}+1,N\rho_i[\, \cap\, \Int)
\end{equation*}
(où l'intervalle $[a,b[$ est vide si $a \soe b$).
\end{prop}
\dem

D'une part, chaque intervalle $[\rho_{i-1}+\delta,\rho_i[$ est inclus dans l'ensemble $X(\delta)$. D'autre part, si~$t \in X(\delta)$, on a~$t<\rho_M=1$ ; soit $i_0$ le plus petit indice $i$ tel que $t< \rho_i$. On a $i_0>0$ et $\rho_{i_0-1} +\delta \ioe t$ (sinon on aurait $\rho_{i_0-1} \in \, ]t-\delta,t]$). Par conséquent $t \in [\rho_{i_0-1}+\delta,\rho_{i_0}[$ et la première égalité ensembliste est démontrée.

Pour la seconde, on observe que
\begin{equation*}
j \in\set{1,\dots,N} \text{ et }  \Big] \frac{j-1}{N}, \frac jN\Big] \, \cap \, F= \vide \; \ssi \; j \in NX(1/N)\, \cap\, \Int.\fine
\end{equation*}

\begin{prop}\label{230131a}
Sous les hypothèses de la proposition \ref{230130a}, on a
\begin{multline*}
\card\Big\{j \in \set{1,\dots,N}  \pepv  \Big] \frac{j-1}{N}, \frac jN\Big] \, \cap \, F= \vide\Big\} =\\
N\mes \big(X(1/N)\big)+\sum_{\substack{i=1\\ \rho_i-\rho_{i-1} \soe 1/N}}^M\,(\set{-N\rho_i} -\set{-N\rho_{i-1}}).
\end{multline*}
\end{prop}
\dem

En notant que les intervalles $[N\rho_{i-1}+1,N\rho_i[$ de la proposition \ref{230130a} (dont certains peuvent être vides) sont deux à deux disjoints, on a
\begin{multline*}
\card\Big\{j \in \set{1,\dots,N}  \pepv  \Big] \frac{j-1}{N}, \frac jN\Big] \, \cap \, F= \vide\Big\} =\\
\sum_{i=1}^M \card([N\rho_{i-1}+1,N\rho_i[\, \cap\, \Int)=\sum_{\substack{i=1\\ \rho_i-\rho_{i-1} \soe 1/N}}^M\, (\lceil N\rho_i \rceil -\lceil N\rho_{i-1}\rceil - 1)\expli{d'après \eqref{230201a}}\\
=N\sum_{\substack{i=1\\ \rho_i-\rho_{i-1} \soe 1/N}}^M\, (\rho_i -\rho_{i-1} - 1/N)+\sum_{\substack{i=1\\ \rho_i-\rho_{i-1} \soe 1/N}}^M\,(\set{-N\rho_i} -\set{-N\rho_{i-1}})\\
=N\mes \big(X(1/N)\big)+\sum_{\substack{i=1\\ \rho_i-\rho_{i-1} \soe 1/N}}^M\,(\set{-N\rho_i} -\set{-N\rho_{i-1}}).\fine
\end{multline*}

\subsection{Identité pour $R(N)$}

Pour $k \soe 1$, nous allons choisir $F=\Fcal_k$ dans la proposition \ref{230131a} ; indiquons la dépendance en $k$ en notant
\[
\rho_0(k)=0 < \rho_1(k) < \dots < \rho_M(k)=1 \quad (M=A(k)),
\]
les éléments de $\Fcal_k \, \cap \, [0,1]$.

\begin{prop}\label{230217b}
Pour $N \in \Nat^*$, on a
\begin{equation*}
R(N)=-2\sum_{k=1}^N\sum_{\substack{i=1\\ \rho_i(k)-\rho_{i-1}(k) \soe 1/N > \rho_{i+1}(k)-\rho_{i}(k)}}^{A(k)-1}\,B_1\big(N\rho_i(k)\big).
\end{equation*}
\end{prop}
\dem

Pour $k=1,\dots, N$, on a
\begin{align*}
\theta_N(k) &=\card\set{j \in\set{1,\dots,N}  \pepv q_j(N) >k}\\
&=\card\Big\{j \in \set{1,\dots,N} \pepv \Big] \frac{j-1}{N}, \frac jN\Big] \, \cap \, \Fcal_k= \vide\Big\} \expli{d'après \eqref{230131b}}\\
&=N\mes \Big(\Big\{t \in\,  ]0,1] \pepv \Big]t- \frac 1N,t\Big] \, \cap \, \Fcal_k= \vide\Big\} \Big)+\\
&\qquad\qquad\sum_{\substack{i=1\\ \rho_i(k)-\rho_{i-1}(k)\soe 1/N}}^{A(k)}\,(\set{-N\rho_i(k)} -\set{-N\rho_{i-1}(k)}),
\end{align*}
d'après la proposition \ref{230131a}.
Or on a
\begin{equation*}
\mes \Big(\Big\{t \in\,  ]0,1] \pepv \Big]t- \frac 1N,t\Big] \, \cap \, \Fcal_k= \vide\Big\}=\mes \big(\set{t \in\,  ]0,1] \pepv q(t,1/N)>k} \big)=\nu_N(k),
\end{equation*}
donc
\[
\theta_N(k)-N\nu_N(k)=\sum_{\substack{i=1\\ \rho_i(k)-\rho_{i-1}(k) \soe 1/N}}^{A(k)}\,(\set{-N\rho_i(k)} -\set{-N\rho_{i-1}(k)}).
\]
En notant $\xi_N(k)$ cette somme, on en déduit avec \eqref{230217a} que
\[
R(N)=\sum_{k=1}^N \xi_N(k).
\]

Pour la fin du raisonnement, le nombre $k$ est fixé et nous écrivons $\rho_i$ au lieu de $\rho_i(k)$, et $M$ au lieu de $A(k)$. La quantité $\xi_N(k)$ se prête à une transformation d'Abel. En posant 
\[
\eps_i=[\rho_i-\rho_{i-1} \soe 1/N] \;\text{ et } \;f(t) =\set{-Nt},
\]
on a
\begin{displaymath}
\xi_N(k)=\sum_{i=1}^{M}\,\big(f(\rho_i)-f(\rho_{i-1})\big)\eps_i
=\sum_{i=1}^{M-1}\,f(\rho_i)(\eps_i-\eps_{i+1}).
\end{displaymath}
Or,
\[
\eps_i-\eps_{i+1}=
\begin{dcases}
0 & \text{ si } \max( \rho_i-\rho_{i-1} , \rho_{i+1}-\rho_{i} ) < 1/N\\
0 & \text{ si } \min( \rho_i-\rho_{i-1} , \rho_{i+1}-\rho_{i} )  \soe 1/N\\
1 &  \text{ si }  \rho_i-\rho_{i-1} \soe 1/N > \rho_{i+1}-\rho_{i} \\
-1&  \text{ si }  \rho_i-\rho_{i-1} < 1/N \ioe  \rho_{i+1}-\rho_{i}.
\end{dcases}
\]
Ainsi,
\begin{equation*}
\xi_N(k)=\sum_{\substack{i=1\\ \rho_i-\rho_{i-1} \soe 1/N > \rho_{i+1}-\rho_{i}}}^{M-1}\,f(\rho_i)\quad -\;\sum_{\substack{i=1\\ \rho_i-\rho_{i-1} < 1/N \ioe  \rho_{i+1}-\rho_{i}}}^{M-1}\,f(\rho_i).
\end{equation*}
Si dans la dernière somme on fait le changement de variable de sommation $j=M-i$, on obtient la somme
\[
\sum_{\substack{j=1\\ \rho_{M-j}-\rho_{M-j-1} < 1/N \ioe  \rho_{M-j+1}-\rho_{M-j}}}^{M-1}\,f(\rho_{M-j})=\sum_{\substack{j=1\\ \rho_{j}-\rho_{j-1} \soe 1/N > \rho_{j+1}-\rho_{j}}}^{M-1}\,f(1-\rho_j),
\]
où l'on a utilisé le fait que $\rho_{M-j}=1-\rho_j$. Par conséquent,
\begin{equation*}
\xi_N(k)=\sum_{\substack{i=1\\ \rho_i-\rho_{i-1} \soe 1/N > \rho_{i+1}-\rho_{i}}}^{M-1}\,\big(f(\rho_i)-f(1-\rho_i)\big).
\end{equation*}
Comme $f(t)=\set{-Nt}$, on a
\[
f(t)-f(1-t)=\set{-Nt} - \set{-N(1-t)}=\set{-Nt} - \set{Nt}=-2B_1(Nt).
\]
Par conséquent,
\begin{equation*}\label{230212a}
\xi_N(k)=-2\sum_{\substack{i=1\\ \rho_i-\rho_{i-1} \soe 1/N > \rho_{i+1}-\rho_{i}}}^{M-1}\,B_1(N\rho_i).\fine
\end{equation*}

\section{Estimation de $R(N)$}\label{sec:esti-RN}

\subsection{Transformation de la somme $R(N)$}

Pour $k,N$ des nombres entiers tels que $1 \ioe k \ioe N$, posons
\begin{equation*}
\sigma_N(k)=\sum_{\substack{i=1\\ \rho_i(k)-\rho_{i-1}(k) \soe 1/N > \rho_{i+1}(k)-\rho_{i}(k)}}^{A(k)-1}\,B_1\big(N\rho_i(k)\big),
\end{equation*}
et
\begin{equation}
  \label{eq:def-TN}
  T(N)=-\frac{R(N)}2=\sum_{k=1}^N \sigma_N(k).
\end{equation}

\begin{prop}\label{prop:decompo-T}
On a
\begin{equation}
  \label{eq:decompo-T}
  T(N)=T(N)=T_1(N)+T_2(N),
\end{equation}
avec
\begin{align*}
T_1(N) &= \sum_{\substack{r<s \\rs \ioe N \\ (r,s)=1}} \; B_1\Big( \frac{N\bar{r}}s\Big)\sum_{\substack{s \ioe k <r+s\\ \thet(k,r,s) >N/s}} 1,\\
T_2(N) &= \sum_{\substack{s<r \\ rs \ioe N \\ (r,s)=1}} \; B_1\Big( \frac{N\bar{r}}s\Big)\sum_{\substack{r \ioe k <r+s\\ \thet(k,r,s) >N/s}} 1.
\end{align*}
\end{prop}
\dem

En utilisant les propriétés rappelées au \S\ref{230131d}, on obtient 
\begin{equation*}
\sigma_N(k) = \sum_{\substack{\max(r,s) \ioe k \\ r+s>k  \\ (r,s)=1 \\  rs\ioe N <s\thet(k,r,s)}}B_1\Big( \frac{N\bar{r}}s\Big),
\end{equation*}
où $\bar{r}=\inv(r,s)$ (observons que le couple $(r,s)=(k,1)$, qui figure dans l'ensemble \eqref{230221b}, ne contribue pas à la somme).
Par conséquent,
\begin{displaymath}
T(N) = \sum_{k=1}^N \;  \sum_{\substack{\max(r,s) \ioe k \\ r+s>k  \\ (r,s)=1 \\  rs\ioe N <s\thet(k,r,s)}}B_1\Big( \frac{N\bar{r}}s\Big)
= \sum_{\substack{rs \ioe N \\ (r,s)=1}} \; B_1\Big( \frac{N\bar{r}}s\Big) \sum_{\substack{k=1\\ \max(r,s) \ioe k <r+s\\\thet(k,r,s) >N/s}}^N 1.
\end{displaymath}
La diagonale $r=s$ ne contribue pas à la somme. On a aussi
\[
rs \ioe N \; \implique \; r+s \ioe r+N/r \ioe N+1,
\]
donc, dans la dernière somme, la condition $k <r+s$  rend la condition $k\ioe N$ superflue.
En séparant la somme suivant les conditions $r<s$ ou $s>r$, on obtient l'énoncé.\fin

\subsection{Décomposition de $T_1(N)$}

Comme 
\[
\thet(k,r,s)=s\Big\lfloor \frac{k+r}s\Big\rfloor-r,
\]
il est naturel de regrouper les termes suivant les valeurs possibles de $j=\lfloor (k+r)/s\rfloor$.

On a
\[
(r <s ) \; \text{ et } \; (s\ioe k < r+s) \; \implique \;  1 + \frac 1s \ioe \frac {k+r}s \ioe 3 - \frac 2s\virg
\]
donc $j$ ne peut prendre que les valeurs $1$ et $2$.
La valeur $j=1$ correspond aux valeurs de~$k$ telles que $s-r \ioe k <2s-r$. La valeur $j=2$ correspond aux valeurs de $k$ telles que~$2s-r \ioe k <3s-r$. Ainsi,
\begin{equation}
  \label{eq:decompo-T1}
  T_1(N)=T_{11}(N)+T_{12}(N),
\end{equation}
avec
\begin{align*}
T_{11}(N)=\sum_{\substack{r<s \\rs \ioe N\\s(s-r)>N \\ (r,s)=1}} \; B_1\Big( \frac{N\bar{r}}s\Big)\sum_{\substack{s \ioe k <\min(r+s,2s-r)}} 1,\\
T_{12}(N)=\sum_{\substack{r<s \\rs \ioe N\\s(2s-r)>N \\ (r,s)=1}} \; B_1\Big( \frac{N\bar{r}}s\Big)\sum_{\substack{2s-r \ioe k <\min(r+s,3s-r)}} 1.
\end{align*}

\subsection{Estimation de $T_{11}(N)$}

En distinguant les valeurs de $r$ et $s$ telles que $r+s \ioe 2s-r$ ou $r+s >2s-r$, on a
\[
T_{11}(N)= \sum_{\substack{r\ioe s/2 \\rs \ioe N\\s(s-r)>N \\ (r,s)=1}} \; B_1\Big( \frac{N\bar{r}}s\Big)\sum_{\substack{s \ioe k <r+s}} 1
+ \sum_{\substack{s/2 <r<s \\rs \ioe N\\s(s-r)>N \\ (r,s)=1}} \; B_1\Big( \frac{N\bar{r}}s\Big)\sum_{\substack{s \ioe k <2s-r}} 1.
\]
Remarquons que 
\begin{displaymath}
  \big(s(s-r)>N\big ) \text{ et } (rs \ioe N)  \Rightarrow  r <\frac 12 \Big(\frac  Ns + s- \frac Ns\Big)=\frac s2
\end{displaymath}
donc la seconde double somme sur $r$ et $s$ est vide, et la condition $r \ioe s/2$ est superflue dans la première.
De plus, comme
\begin{displaymath}
  (s(s-r)>N) \text{ et } (r >0)  \Rightarrow s>\sqrt{N}, 
\end{displaymath}
on a
\begin{align*}
T_{11}(N)
& =\sum_{s >  \sqrt{N}} \;\sum_{\substack{1\le r \ioe \min(N/s,s-(N+1)/s) \\ (r,s)=1}} \; rB_1\Big( \frac{N\bar{r}}s\Big). 
\end{align*}
D'après la proposition \ref{eval-somme-B1-IPP}, on a
\begin{displaymath}
  \Big|\sum_{\substack{1\le r \ioe \min(N/s,s-(N+1)/s) \\ (r,s)=1}} \; rB_1\Big( \frac{N\bar{r}}s\Big)\Big| \le \frac N s   \sqrt{s}\, \tau(s)\beta(s)\ln s
= \frac{N}{\sqrt{s}} \tau(s)\beta(s)\ln s. 
\end{displaymath}
Soit $K$ un nombre réel tel que $\sqrt{N}\, \ioe K \ioe N$. Pour $ K< s \le N$, nous utilisons la majoration triviale 
\begin{displaymath}
  \Big|\sum_{\substack{1\le r \ioe \min(N/s,s-(N+1)/s) \\ (r,s)=1}} \; rB_1\Big( \frac{N\bar{r}}s\Big)\Big| \le 
\sum_{1\le r\ioe N/s} r \le \frac{N^2}{s^2}. 
\end{displaymath}
Nous avons donc 
\begin{align*}
 | T_{11}(N)|  \le  \sum_{\sqrt{N} < s \le  K} \frac{N}{\sqrt{s}}  \tau(s)\beta(s)\ln s 
+ N^2\sum_{ K< s \le N} \frac{1}{s^2}.  
\end{align*}
L'estimation \eqref{230222a} permet d'obtenir, par sommation partielle, la majoration
\begin{displaymath}
  \sum_{n\le t} \frac{\tau(n)\beta(n)\ln n}{\sqrt{n}} =O(\sqrt{t} \ln^3 t)  \quad(t\ge 1). 
\end{displaymath}
Nous obtenons ainsi
\[
  T_{11}(N) \ll NK^{1/2}\ln^3 N + N^2/K.
\]
Le choix $K=N^{2/3}/\ln^2 N$ (pour $N$ assez grand) fournit alors l'estimation 
\begin{equation}\label{eq:majo-T11}
T_{11}(N) \ll N^{4/3}\ln^2 N.
\end{equation}

\subsection{Estimation de $T_{12}(N)$}
Comme 
\begin{displaymath}
  r<s \Rightarrow r+s<3s-r,
\end{displaymath}
on a 
\begin{displaymath}
  T_{12}(N) =\sum_{\substack{r<s \\rs \ioe N\\s(2s-r)>N \\ (r,s)=1}} \; B_1\Big( \frac{N\bar{r}}s\Big)\sum_{2s-r \ioe k <r+s} 1  = \sum_{\substack{s/2<r \ioe s \\rs \ioe N\\s(2s-r)>N \\ (r,s)=1}} \; (2r-s) B_1\Big( \frac{N\bar{r}}s\Big).
\end{displaymath}
Notons que 
\begin{displaymath}
  (rs\le N) \text{ et } ( s/2<r ) \Rightarrow s \ioe \sqrt{2N}. 
\end{displaymath}
Par conséquent, 
\begin{displaymath}
  T_{12}(N) =2
\sum_{1\le s \le \sqrt{2N}}\,\sum_{\substack{s/2 < r \ioe \min(s,N/s,2s -(N+1)/s)\\(r,s)=1}} \; (r-s/2) B_1\Big( \frac{N\bar{r}}s\Big).
\end{displaymath}
D'après la proposition \ref{eval-somme-B1-IPP}, on a 
\begin{displaymath}
\Big| \sum_{s/2 <r \ioe \min(s,N/s, 2s-(N+1)/s)} \; (r-s/2) B_1\Big( \frac{N\bar{r}}s\Big)\Big|\le  \frac{s^{3/2}}2 \tau(s)\beta(s) \ln s,
\end{displaymath}
et donc 
\begin{align}
  \abs{T_{12}(N)} &\ioe \sum_{1\le s \le \sqrt{2N}}  s^{3/2} \tau(s)\beta(s) \ln s\notag\\
& \ll N^{3/4} \ln N  \sum_{1\le s \le \sqrt{2N}} \tau(s)\beta(s)\notag\\
& \ll N^{5/4} \ln^3 N \expli{d'après \eqref{230222a}}. \label{eq:majo-T12}
\end{align}

\subsection{Estimation finale de $T_{1}(N)$}

En insérant \eqref{eq:majo-T11} et \eqref{eq:majo-T12} dans \eqref{eq:decompo-T1}, nous obtenons 
\begin{equation}\label{eq:majo-T1}
  T_1(N) \ll N^{4/3} \ln^2 N.  
\end{equation}

\subsection{\'Evaluation de $T_{2}(N)$}
Rappelons que
\begin{displaymath}
T_2(N) = \sum_{\substack{s<r \\ rs \ioe N \\ (r,s)=1}} \; B_1\Big( \frac{N\bar{r}}s\Big)\sum_{\substack{r \ioe k <r+s\\ \thet(k,r,s) >N/s}} 1.
\end{displaymath}
Nous regroupons  les termes suivant les valeurs possibles de $j=\lfloor (k+r)/s\rfloor$.
On a
\[
(s<r ) \; \text{ et } \; (s\ioe k < r+s) \; \implique \;  2 + \frac 1s \ioe \frac {k+r}s < \frac {2r}s +1 \le 2N+1.
\]
donc $j$ peut prendre les valeurs $2,\ldots, 2N$.
Comme 
\begin{displaymath}
  \Big\lfloor \frac{k+r}s\Big\rfloor =j \Leftrightarrow js-r\ioe  k<(j+1)s-r, 
\end{displaymath}
on en déduit
\begin{displaymath} 
T_2(N) = \sum_{j=2}^{2N} \sum_{\substack{s<r \\ rs \ioe N \\ s(js-r)>N\\ (r,s)=1}} \; B_1\Big( \frac{N\bar{r}}s\Big)\sum_{\substack{\max(r, js-r) \ioe k <\min(r+s,(j+1)s-r) }} 1. 
\end{displaymath}
On remarque que 
\begin{displaymath}
  rs \le N \text{ et } js-r> \frac Ns \Rightarrow js-r >r,
\end{displaymath}
donc 
\begin{align} 
  T_2(N) &= \sum_{j=2}^{2N} \sum_{\substack{s<r \\ rs \ioe N \\ s(js-r)>N\\ (r,s)=1}} \; B_1\Big( \frac{N\bar{r}}s\Big)\sum_{\substack{js-r \ioe k <\min(r+s,(j+1)s-r) }} 1\notag \\
         & = \sum_{j=2}^{2N} \sum_{\substack{s<r \\ rs \ioe N \\ r\le js/2\\ s(js-r)>N\\ (r,s)=1}} \; B_1\Big( \frac{N\bar{r}}s\Big)\sum_{\substack{js-r \ioe k <r+s }} 1 \label{T2-1}\\
         & \qquad \qquad
           +  \sum_{j=2}^{2N} \sum_{\substack{s<r \\ rs \ioe N \\ r> js/2\\ s(js-r)>N\\ (r,s)=1}} \; B_1\Big( \frac{N\bar{r}}s\Big)\sum_{\substack{js-r \ioe k <(j+1)s-r }} 1.\label{T2-2}
\end{align}
La double somme intérieure en $r$ et $s$ dans \eqref{T2-2} est vide car les conditions sur $r$ et $s$ entraînent 
\begin{displaymath}
  N < s(js-r) <s (js-js/2) = \frac{js}{2} s \le  \frac{js}{2} \frac Nr < N.  
\end{displaymath}
De plus, pour que la somme intérieure soit non vide dans \eqref{T2-1}, on doit avoir $js-r<r+s$, \cad~$ r >(j-1)s/2$. Par conséquent, 
\begin{displaymath}
  T_{2}(N) =  \sum_{j=2}^{2N} \sum_{\substack{s<r \\ rs \ioe N \\ (j-1)s/2< r\le js/2\\ s(js-r)>N\\ (r,s)=1}} \; \big(2r-(j-1)s\big)B_1\Big( \frac{N\bar{r}}s\Big).
\end{displaymath}
Le terme correspondant à $j=2$ vaut $0$ puisque les conditions entraînent $s<r \le s$. 
Par ailleurs, 
\begin{displaymath}
   \frac{(j-1)s}2< r \text{ et } r\le N/s \Rightarrow s \le \sqrt{\frac{2N}{j-1}}.
\end{displaymath}
Ainsi 
\begin{displaymath}
  T_{2}(N) = 2 \sum_{j=3}^{2N} \sum_{ s \le  \sqrt{2N/(j-1)}}\;\; \sum_{\substack{(j-1)s/2<r \le \min(js/2, N/s, js - (N+1)/s ) \\ (r,s)=1}} \; \big(r-(j-1)s/2\big)B_1\Big( \frac{N\bar{r}}s\Big).
\end{displaymath}
On constate que la sommation en $r$ porte sur un intervalle (peut-être vide) dont la longueur n'excède pas $s/2$. 
On a donc d'après la proposition \ref{eval-somme-B1-IPP}
\begin{displaymath}
 \Big| \sum_{\substack{(j-1)s/2<r \le \min(js/2, N/s, js - (N+1)/s ) \\ (r,s)=1}} \; \big(r-(j-1)s/2)\big)B_1\Big( \frac{N\bar{r}}s\Big)\Big|
\le \frac{s}2 \sqrt{s}\, \tau(s)\beta(s) \ln s. 
\end{displaymath}
Il suit 
\begin{align}
  T_{2}(N) &\le \sum_{j=3}^{2N} \,\sum_{1 \le s \le  \sqrt{2N/(j-1)}}  s^{3/2} \tau(s)\beta(s) \ln s\notag 
\\& \ll  \ln^3 N  \sum_{j=3}^{2N}   \Big(\frac{N}{j-1}\Big)^{5/4} \expli{cf. \eqref{eq:majo-T12}}\notag\\& \ll N^{5/4} \ln^3 N.\label{eq:majo-T2}
\end{align}

\subsection{Conclusion}
En insérant \eqref{eq:majo-T1} et \eqref{eq:majo-T2} dans \eqref{eq:decompo-T}, nous obtenons l'estimation 
\begin{displaymath}
  T(N) \ll N^{4/3}  \ln^2 N
\end{displaymath}
et donc, d'après  \eqref{eq:def-TN},
\begin{displaymath}
  R(N) \ll N^{4/3}  \ln^2 N.
\end{displaymath}
Cela conclut la preuve du théorème. 

\section{Variantes du théorème}\label{sec:variante-th}

Posons
\begin{align*}
\bar{S} (N) &= \sum_{j=1}^N q\Big(\,\Big[\frac{j-1}N, \frac jN \Big]\,\Big) ,\\
S^* (N) &= \sum_{j=1}^N q\Big(\,\Big[ \frac{j-1}N, \frac jN \Big[\,\Big) ,\\
\tilde{S}(N)&= \sum_{j=1}^N q\Big(\,\Big] \frac{j-1}N, \frac jN \Big[\,\Big) .
\end{align*}
Nous démontrons dans cette section que la formule asymptotique   \eqref{eq:resultat-principal} reste valable pour les quantités $\bar{S} (N)$, $S^* (N)$ et $\tilde{S}(N)$. 

D'abord, la fonction $x\mapsto 1-x$ transforme $[(j-1)/N,j/N[$ en $](N-j)/N,(N-j+1/N]$ et conserve les dénominateurs minimaux. Par conséquent,
\[
q\Big( \Big[ \frac{j-1}{N}, \frac jN\Big[\Big)=q\Big(\Big]\frac{N-j}N, \frac{N-j+1}N\Big]\Big) \quad (1 \ioe j \ioe N),
\]
et donc $S^*(N)=S(N)$.
Ensuite, comme $E \mapsto q(E)$ est décroissante, on a
\[
\bar{S} (N)  \ioe S^* (N) =S(N) \ioe \tilde{S}(N).
\]
\begin{prop}
  On a pour $N\ge 1$, 
  \begin{displaymath}
     \tilde{S}(N)-S(N) \le N \tau(N) \quad \text{ et } \quad   S(N)- \bar{S}(N)\le N \tau(N). 
  \end{displaymath}
\end{prop}
Comme $\tau(N) \ll_\eps N^\eps$, cela entraîne bien que $\bar{S}(N)$ et $\tilde{S}(N)$ satisfont également \eqref{eq:resultat-principal}. 
\medskip

\dem 
Nous démontrons seulement la première inégalité, la seconde se démontre de même. 
Introduisons les notations  
\begin{displaymath}
  \tilde{q}_j(N) = q\Big(\Big] \frac{j-1}N, \frac jN \Big[\Big) 
\quad \text{ et } \quad
\tilde{\theta}_N(k) = \card\{ j\in\{1,\ldots,N\} \pepv   \tilde{q_j}(N) >k\}.
\end{displaymath}

Notons quelques différences entre ces quantités et $q_j(N)$ et $\theta_N(k)$. On n'a jamais la relation~$\tilde{q}_j(N)=N$ ; en revanche, on a $\tilde{q}_1(N)=\tilde{q}_N(N)=N+1$ et
\[
\tilde{\theta}_N(N-1) =\tilde{\theta}_N(N) =2 \epv \tilde{\theta}_N(N+1) =0 \quad (N \soe 2).
\]

Néanmoins, en procédant comme pour \eqref{eq:identite-SN}, nous obtenons l'identité 
\begin{displaymath}
   \tilde{S}(N)=\sum_{k=0}^N\, \tilde{\theta}_N(k), 
\end{displaymath}
d'où
\begin{displaymath}
   \tilde{S}(N)-S(N)=\sum_{k=1}^N \big(\tilde{\theta}_N(k)-{\theta}_N(k)\big).
\end{displaymath}
Comme $q_j(N)\le \tilde{q_j}(N)$, on a 
\begin{displaymath}
   \tilde{S}(N)-S(N)
= \sum_{k=1}^N
\card\{ j\in\{1,\ldots,N\} \pepv    \tilde{q_j}(N)>k \,\text{ et }  \,q_j(N)\le k\}.
\end{displaymath}
Soit $j\in\{1,\ldots,N\}$ tel que $ \tilde{q_j}(N)>k \text{ et } q_j(N)\le k$. Il existe donc $i\in\{1,\ldots,A(k)\}$ tel que 
\begin{displaymath}
 \frac{j}{N} = \rho_i(k)  \quad \text{ et } \quad \rho_{i-1}(k)\le \frac{j-1}{N}, 
\end{displaymath}
ce qui implique 
\begin{equation}\label{eq:carac-ij}
  N\rho_{i}(k) \in\Int  \quad \text{ et }  \quad \rho_{i}(k)-\rho_{i-1}(k)\ge \frac 1N.
\end{equation}
Réciproquement, si $i\in\{1,\ldots,A(k)\}$ est tel que \eqref{eq:carac-ij}, il existe $j\in\{1,\ldots,N\}$ tel que~$ \tilde{q_j}(N)>k \text{ et } q_j(N)\le k$. 
On a donc établi que
\begin{displaymath}
   \tilde{S}(N)-S(N)
= \sum_{k=1}^N \sum_{\substack{i=1\\  \rho_{i}(k)-\rho_{i-1}(k)\ge \frac 1N }}^{A(k)} [N \rho_i(k)\in\Int].  
\end{displaymath}
En utilisant la paramétrisation de $\Fcal_k\, \cap \, ]0,1]$ donnée par \eqref{230221b}, on obtient
\begin{align*}
\tilde{S} (N)&-S(N) =\sum_{k=1}^N\,\sum_{\substack{\max(r,s) \ioe k\\ r+s >k\\rs \ioe N\\(r,s)=1}}[s \mid N]\\
&=\sum_{\substack{rs \ioe N\\(r,s)=1}}([s \mid N])\sum_{ \max(r,s) \ioe k < r+s}1\expli{on a $r+s \ioe N+1$, cf. preuve de la proposition~\ref{prop:decompo-T}}\\
&=\sum_{\substack{rs \ioe N\\(r,s)=1}}[s \mid N])\min(r,s).
\end{align*}
Or
\begin{displaymath}
\sum_{\substack{rs \ioe N\\(r,s)=1}}[s \mid N] \min(r,s) = 
\sum_{s \mid N} \sum_{\substack{r \ioe N/s\\(r,s)=1}} \min(r,s)
\ioe \sum_{s \mid N} s\cdot N/s= N\tau(N),
\end{displaymath}
et donc 
\begin{displaymath}
   \tilde{S}(N)-S(N) \le N\tau(N). \fine
\end{displaymath}

\begin{center}
{\sc Remerciements }
\end{center}
Pour la préparation de ce travail, les deux auteurs ont bénéficié du programme   Research in Paris à l'Institut Henri Poincaré. Ils remercient chaleureusement cette institution pour l'excellent accueil qu'ils ont reçu.


\providecommand{\bysame}{\leavevmode ---\ }
\providecommand{\og}{``}
\providecommand{\fg}{''}
\providecommand{\smfandname}{et}
\providecommand{\smfedsname}{\'eds.}
\providecommand{\smfedname}{\'ed.}
\providecommand{\smfmastersthesisname}{M\'emoire}
\providecommand{\smfphdthesisname}{Th\`ese}

\medskip

\footnotesize

\noindent BALAZARD, Michel\\
Aix Marseille Univ, CNRS, I2M, Marseille, France\\
Adresse \'electronique : \texttt{balazard@math.cnrs.fr}

\medskip

\noindent MARTIN, Bruno\\
ULCO,  LMPA, Calais, France\\
Adresse \'electronique : \texttt{Bruno.Martin@univ-littoral.fr}

\end{document}